\documentclass[a4paper,12pt]{article}
\usepackage{amsmath}
\usepackage{amsthm}
\usepackage{amssymb}
\usepackage{amscd}

\title{The Higher Spin Dirac Operators on 3-Dimensional Manifolds}
\author{Yasushi Homma \thanks{Department of Mathematics, 
  Waseda University, 
  Tokyo 169, 
  Japan, \endgraf 
  {\it E-mail address}: homma{\char'100}gm.math.waseda.ac.jp}}
\date{}

\theoremstyle{plain}
\newtheorem{thm}{Theorem}[section]
\newtheorem{prop}[thm]{Proposition}
\newtheorem{cor}[thm]{Corollary}
\newtheorem{lem}[thm]{Lemma}
\theoremstyle{definition}

\numberwithin{equation}{section}

\theoremstyle{remark}
\newtheorem{rem}{Remark}[section]
\newtheorem{ex}{Example}[section]

\allowdisplaybreaks[3]
\begin{document}
\maketitle
\begin{abstract}
We study the higher spin Dirac operators on $3$-dimensional manifolds and show that there exist two Laplace type operators for each associated bundle. Furthermore, we give lower bound estimations for the first eigenvalues of these Laplace type operators.
\end{abstract}
\section{Introduction}\label{sec:0}
In this paper, we study the higher spin Dirac operator, which is a generalization of the Dirac operator as follows (see \cite{F}, \cite{BS}, and \cite{Bu}). Let $M$ be a $n$-dimensional spin manifold and $\mathbf{Spin}(M)$ be the principal spin bundle on $M$. The irreducible unitary representation $(\rho, V_{\rho})$ of the structure group $Spin(n)$ induces the associated (irreducible) bundle $\mathbf{S}_{\rho}(M)$, 
\begin{equation} 
\mathbf{S}_{\rho}(M):=\mathbf{Spin}(M)\times_{\rho}V_{\rho}.
\end{equation}
For each bundle, we have the covariant derivative $\nabla$ associated to the Levi-Civita connection or the spin connection,
\begin{equation}
\nabla:\Gamma (\mathbf{S}_{\rho}(M))\to \Gamma (\mathbf{S}_{\rho}(M)\otimes T^{\ast}(M)).
\end{equation}
Here, the cotangent bundle $T^{\ast}(M)$ is the bundle corresponding to the adjoint representation $(\mathrm{Ad},\mathbf{R}^n)$ of $Spin(n)$. So we decompose the tensor bundle $\mathbf{S}_{\rho}(M)\otimes T^{\ast}(M)$ into irreducible bundles with respect to $Spin(n)$. Let $\pi_{\rho, \nu}$ be the orthogonal projection to the irreducible bundle $\mathbf{S}_{\nu}(M)$ from $\mathbf{S}_{\rho}(M)\otimes T^{\ast}(M)\simeq \oplus_{\nu}\mathbf{S}_{\nu}(M)$. Then we define the higher spin Dirac operator $D_{\rho, \nu}$ to be the composed mapping $\pi_{\rho, \nu}\circ \nabla$,
\begin{equation}
D_{\rho, \nu}:\Gamma (\mathbf{S}_{\rho}(M))\xrightarrow{\nabla} \Gamma (\mathbf{S}_{\rho}(M)\otimes T^{\ast}(M))\xrightarrow{\pi_{\rho, \nu}}\Gamma (\mathbf{S}_{\nu}(M)).
\end{equation}
In fact, the Dirac operator is given in this way. To construct the Dirac operator, we take the spinor representation $(\Delta, V_{\Delta})$ and the associated bundle $\mathbf{S}_{\Delta}(M)$. Then the tensor bundle $\mathbf{S}_{\Delta}(M)\otimes T^{\ast}(M)$ decomposes into the direct sum of only two irreducible bundles, $\mathbf{S}_{\Delta}(M)$ and $\mathbf{S}_{T}(M)$. Then the differential operator $D:=D_{\Delta,\Delta}$ is the Dirac operator and $D_{\Delta,T}$ is the twistor operator (see \cite{BFGK} and \cite{BF}). On the other hand, we know another definition of the Dirac operator by using the Clifford algebra, that is,
\begin{equation}
D=\sum_i e_i\cdot\nabla_{e_i}.\label{eqn:0-3}
\end{equation}
From the relations 
\begin{equation}
e_i e_j+e_j e_i=-2\delta_{ij}, \label{eqn:0-2}
\end{equation}
we show that the Dirac operator satisfies the Bochner type identity 
\begin{equation}
D^2=\nabla^{\ast}\nabla+\frac{1}{4}\kappa, \label{eqn:0-1}
\end{equation}
where $\kappa$ is the scalar curvature of $M$. 

The aim of this paper is to give the Bochner type identities for the higher spin Dirac operators on $3$-dimensional spin manifolds. As mentioned above, the relations \eqref{eqn:0-2} is necessary to give the Bochner type identity for the Dirac operator. But the Clifford action does not exist on the representation spaces of $Spin(n)$ in general. So we consider linear mappings among the representation spaces, which are called {\it the Clifford homomorphisms}. For a $3$-dimensional spin manifold $M$, the structure group of $\mathbf{Spin}(M)$ is $Spin(3)=SU(2)$. Then we use the Clebsch-Gordan formula to define the Clifford homomorphisms. By using the Clifford homomorphisms, we obtain local formulas of the higher spin Dirac operators such as \eqref{eqn:0-3} and the Bochner type identities for them. Furthermore, the identities lead us to give lower bound estimations for the first eigenvalues of these operators.

In section \ref{sec:1}, we explain the Clebsch-Gordan formula for the Lie group $SU(2)$. In section \ref{sec:2}, we define the Clifford homomorphisms on the representation spaces and obtain some relations among these homomorphisms including the usual Clifford relations \eqref{eqn:0-2}. In section \ref{sec:3}, we have formulas of the higher spin Dirac operators by using the Clifford homomorphisms and investigate the properties of these operators (ellipticity, the Bochner type identities, and so on.). The interest thing is that we obtain two Laplace type operators for each associated bundle. In section \ref{sec:4}, we have the lower bound estimations for the first eigenvalue of the Laplace type operators. This estimation is a generalization of the one for the Dirac operator given in \cite{BFGK} or the Laplace-Beltrami operator in \cite{GM} and \cite{P}. In the section \ref{sec:5}, we consider the case of the $3$-dimensional manifold of the constant curvature and show that some operators commute. In the last section, as an example, we calculate all the eigenvalues of the higher spin Dirac operators on the symmetric space $S^3$.  
\section{The Clebsch-Gordan formula }\label{sec:1}
In this section we shall explain the representations of $SU(2)$ and the Clebsch-Gordan formula. Let $V_m$ be the $(m+1)$-dimensional complex vector space of polynomials of degree $\le m$ in $z_m$. The inner product on $V_m$ is set by
\begin{equation}
        \left(v_m^k,v_m^l\right)=\delta_{kl}, \label{1-1}
\end{equation}
where 
\begin{equation}
v_m^k:=\frac{z^k_m}{\sqrt{k!(m-k)!}}.  \label{1-2}
\end{equation}
We define a representation $\rho_{m}$ on $V_m$ by $\rho_{m}(h)z_m^k=(b z_m+d)^{m-k}(a z_m+c)^k$ for $h=
\left( 
\begin{smallmatrix}
a & b \\
c & d
\end{smallmatrix}
\right) 
$ in $SU(2)$.
Then $(\rho_m, V_m)$ is a finite dimensional irreducible unitary representation of $SU(2)$ called the spin-$\frac{m}{2}$ representation and all such representations are given in this way. 

We denote the infinitesimal representation of $(\rho_m, V_m)$ by the same symbol $(\rho_m, V_m)$. The Lie algebra $\mathfrak{su}(2)$ of $SU(2)$ has the following basis, that is, the Pauli matrices:
    \begin{equation}
     \sigma _1:=
         \begin{pmatrix} 
                 i & 0 \\
                 0 & -i 
          \end{pmatrix}, 
   \quad
     \sigma _2:=
         \begin{pmatrix} 
               0 & 1 \\ 
               -1 &  0 
       \end{pmatrix},    
 \quad
  \sigma _3:=
       \begin{pmatrix} 
            0 & i \\ 
            i &  0 
       \end{pmatrix}.  \label{1-3}
   \end{equation} 
Then we show that 
\begin{align}
\rho_{m}(\frac{\sigma_1}{2})z_m^k &=i(k-\frac{m}{2})z_m^k, \nonumber \\
\rho_{m}(\frac{\sigma_2}{2}+i\frac{\sigma_3}{2})z_m^k &=-k z_m^{k-1}, \label{1-4} \\
\rho_{m}(\frac{\sigma_2}{2}-i\frac{\sigma_3}{2})z_m^k &=(m-k)z^{k+1}_m. \nonumber
\end{align}
\begin{ex}
The spin-$\frac{1}{2}$ representation $(\rho_1,V_1)$ is the spinor representation on $\mathbf{C}^2$, where we identify $Spin(3)$ with $SU(2)$. 
\end{ex}
\begin{ex}
The spin-$1$ representation $(\rho_2,V_2)$ is the adjoint representation on $\mathfrak{su}(2)\otimes \mathbf{C}$ of $SU(2)$, or the adjoint representation on $\mathbf{R}^3\otimes \mathbf{C}$ of $Spin(3)$. Here, the correspondence of the bases is given as follows: 
\begin{equation}
z_2^0\leftrightarrow\frac{\sigma_2+i\sigma_3}{2}\leftrightarrow e_2+i e_3,
\quad
z^1_2\leftrightarrow \frac{i\sigma_1}{2}\leftrightarrow i e_1,
\quad
z_2^2\leftrightarrow\frac{\sigma_2-i\sigma_3}{2}\leftrightarrow e_2-i e_3, \label{1-5}
\end{equation}
where $z_2^i$ is in $V_2$, $\sigma_i$ in $\mathfrak{su}(2)$, and $e_i$ in $\mathbf{R}^3$.
\end{ex}

Now, we consider the unitary representation $(\rho_m\otimes \rho_n, V_m\otimes V_n)$. Then we can decompose $\rho_m\otimes \rho_n$ into its irreducible components, 
\begin{equation}
\rho_m\otimes \rho_n\simeq \rho_{m+n}\oplus \rho_{m+n-2}\oplus \cdots \oplus \rho_{|m-n|}. \label{1-6}
\end{equation}
This formula is called Clebsch-Gordan formula. We need the orthogonal projection to each irreducible component from $V_m\otimes V_n$ in the next section. 
\section{The Clifford homomorphisms}\label{sec:2}
In this section we shall define the Clifford homomorphisms, which is a generalization of the Clifford action. 
Let $\mathbf{C}l_3$ be the complex Clifford algebra associated to $\mathbf{R}^3$ and $\{e_i\}_{1\le i \le 3}$ be the standard basis of $\mathbf{R}^3$. 
We realize $\mathbf{C}l_3$ as matrix algebra $\mathbf{C}(2)\oplus \mathbf{C}(2)$ by the mapping 
\begin{equation}
\mathbf{C}l_3 \ni e_i\mapsto (\sigma_i,-\sigma_i) \in \mathbf{C}(2)\oplus \mathbf{C}(2).        \label{eqn:2-1}
\end{equation}
Then the Clifford action of $e_i$ on the spinor space $V_1\simeq \mathbf{C}^2$ is given by $e_i\cdot v=\sigma_i v$. Since we would like to generalize this Clifford action on other representation spaces, we use another definition of the Clifford action as follows: we recall the irreducible decomposition
\begin{equation}
(\rho_1, V_1)\otimes (\rho_2, V_2)\simeq (\rho_{3}, V_3)\oplus (\rho_{1}, V_1) \label{eqn:2-2}
\end{equation}
and the isomorphism
\begin{equation}
(\rho_2,V_2)\simeq (ad, \mathfrak{su}(2)\otimes \mathbf{C}) \simeq (ad, \mathbf{R}^3\otimes \mathbf{C}).  \label{eqn:2-3}
\end{equation}
For $v$ in $V_1$ and $e_i$ in $\mathbf{R}^3$, we project $v\otimes e_i$ onto $V_1$ along $V_3$ orthogonally. By calculating the Clebsch-Gordan coefficients, we show that $\mathrm{pr}(v\otimes e_i)=\sigma_i v=e_i\cdot v$.

Now, we consider the representation space $V_m$. In this case, we use the irreducible decomposition
\begin{equation}
(\rho_m, V_m)\otimes (\rho_2, V_2)\simeq (\rho_{m+2}, V_{m+2})\oplus (\rho_m, V_m) \oplus (\rho_{m-2}, V_{m-2}).  \label{eqn:2-4}
\end{equation}
For $v$ in $V_m$ and $X$ in $\mathbf{R}^3$, we decompose $v\otimes X$ as
\begin{equation}
v\otimes X =(v\otimes X)^{+}+(v\otimes X)^{0}+(v\otimes X)^{-}.  \label{eqn:2-5}\end{equation}
Here, $(v\otimes X)^{0}$ is in $V_m$ and $(v\otimes X)^{\pm}$ in $V_{m\pm 2}$. Thus, we have linear mappings from $V_m$ to $V_m$ or $V_{m\pm 2}$ for any $X$ in $\mathbf{R}^3$:
\begin{align}
\rho_m^0(X)v:=&-\frac{\sqrt{m(m+2)}}{2}(v\otimes X)^{0} \in V_{m}, 
                         \nonumber \\
\rho_m^{+}(X)v:=&-\frac{\sqrt{(m+1)(m+2)}}{\sqrt{2}}(v\otimes X)^{+} \in V_{m+2},           \label{eqn:2-6}\\
\rho_m^{-}(X)v:=&\frac{\sqrt{m(m+1)}}{\sqrt{2}}(v\otimes X)^{-} \in V_{m-2},   \nonumber
\end{align}
where we multiply each mapping by a constant to let the calculations easier. We call these linear mappings {\it the Clifford homomorphisms}. 

Calculating the Clebsch-Gordan coefficients in the decomposition \eqref{eqn:2-4}by using {\it Mathematica} (see\cite{Mathe}), we deduce explicit formulas of the Clifford homomorphisms.
\begin{prop}\label{prop:2-1}
The Clifford homomorphisms associated to $\mathbf{R}^3$ are given as follows: for the basis $\{ z_m^k\}_{0\le k \le m}$ of $V_m$ and $\{e_i\}_{1 \le i \le 3}$ in $\mathbf{R}^3$, 
\begin{enumerate}
\item $\rho_m^0(\cdot):V_m\to V_m$,
\begin{align}
\rho_{m}^0(\frac{e_1}{2})z^k_m &=i(k-\frac{m}{2})z_m^k, \nonumber \\
\rho_{m}^0(\frac{e_2}{2}+i\frac{e_3}{2})z^k_m &=-k z^{k-1}_m, \label{eqn:2-7} \\
\rho_{m}^0(\frac{e_2}{2}-i\frac{e_3}{2})z^k_m &=(m-k)z^{k+1}_m. \nonumber
\end{align}
\item 
$\rho_m^{+}(\cdot):V_m\to V_{m+2}$,
\begin{align}
\rho_{m}^+(\frac{e_1}{2})z^k_m &=i z^{k+1}_{m+2}, \nonumber \\
\rho_{m}^+(\frac{e_2}{2}+i\frac{e_3}{2})z^k_m &=-z^{k}_{m+2}, \label{eqn:2-8} \\
\rho_{m}^+(\frac{e_2}{2}-i\frac{e_3}{2})z^k_m &=-z^{k+2}_{m+2}. \nonumber
\end{align}
\item
$\rho_m^{-}(\cdot):V_m\to V_{m-2}$,
\begin{align}
\rho_{m}^-(\frac{e_1}{2})z^k_m &=ik(m-k)z^{k-1}_{m-2}, \nonumber \\
\rho_{m}^-(\frac{e_2}{2}+i\frac{e_3}{2})z^k_m &=k(k-1)z^{k-2}_{m-2}, \label{eqn:2-9} \\
\rho_{m}^-(\frac{e_2}{2}-i\frac{e_3}{2})z^k_m &=(m-k)(m-k-1)z^{k}_{m-2}. \nonumber
\end{align}
\end{enumerate}
\end{prop}
We remark that $\rho_m^0$ is the representation $(\rho_m,V_m)$ of $\mathfrak{su}(2)$ under the isomorphism $\mathfrak{su}(2)\simeq \mathbf{R}^3$ and $\rho_1^0$ is the usual Clifford action on the spinor space $V_1$.

Now, we shall investigate some properties of the Clifford homomorphisms. 
\begin{lem}\label{lem:2-2}
For $X$ in $\mathbf{R}^3\simeq \mathfrak{su}(2)$, we have 
\begin{align}
(\rho_m^0(X))^{\ast}&=-\rho_m^0(X), \label{eqn:2-10}\\
(\rho_m^{\pm}(X))^{\ast}&=-\rho_{m\pm 2}^{\mp}(X), \label{eqn:2-11}
\end{align}
where $(\cdot)^{\ast}$ is the transposed conjugate with respect to the inner product of each $V_m$.
\end{lem}
\begin{proof}
Because $\rho_m^0$ is the representation of $\mathfrak{su}(2)$, the relation \eqref{eqn:2-10} is trivial. So we shall prove that $(\rho_m^+(X))^{\ast}=-\rho_{m+2}^-(X)$. We take the complexification of \eqref{eqn:2-11} and may prove $(\rho_m^+(X+i Y))^{\ast}=-\rho_{m+2}^-(X)+i\rho_{m+2}^-(Y)$. For example, we have
\begin{equation}
\begin{split}
(\rho_{m}^+(\frac{\sigma_2}{2}-i\frac{\sigma_3}{2})z_m^k,z_{m+2}^l) 
&=-(z^{k+2}_{m+2},z^l_{m+2})   \\
&=-(k+2)!(m-k)!\delta_{k+2,l},\quad \mbox{for any $k,l$}.
\end{split}  \nonumber
\end{equation}
On the other hands, 
\begin{equation}
\begin{split}
(z^k_m,-\rho_{m}^-(\frac{\sigma_2}{2}+i\frac{\sigma_3}{2})z^l_{m+2}) 
&=-l(l-1)(z_{m}^{k},v_{m}^{l-2})\\
&=-l(l-1)k!(m-k)!\delta_{k,l-2} \\
&=-(k+2)!(m-k)!\delta_{k+2,l}, \quad \mbox{for any $k,l$}.
\end{split} \nonumber 
\end{equation}
So we have $(\rho_{m}^+(\sigma_2-i\sigma_3))^{\ast}=-\rho_{m}^-(\sigma_2)+i\rho_{m}^-(\sigma_3)$. Similarly we can prove the other cases. 
\end{proof}
\begin{lem}\label{lem:2-3}
For $X$ in $\mathbf{R}^3\simeq \mathfrak{su}(2)$ and $g$ in $SU(2)$, we have 
\begin{align}
\rho_m^0(gXg^{-1})&=\rho_m(g)\rho_m^0(X)\rho_m(g^{-1}) , \label{eqn:2-12}\\
\rho_m^{\pm}(gXg^{-1})&=\rho_{m\pm 2}(g)\rho_m^{\pm}(X)\rho_m(g^{-1}) .\label{eqn:2-13}
\end{align}
\end{lem}
\begin{proof}
 The equation \eqref{eqn:2-12} is trivial. So we shall prove \eqref{eqn:2-13}. For an orthonormal basis $\{v_{m+2}^{k}\}$ of $V_{m+2}$, we denote the corresponding one of the irreducible component $V_{m+2}$ in $V_m\otimes V_2$ by $\{\omega_{m+2}^k\}_k$. Since $\rho_m^{+}$ is the orthogonal projection from $V_m\otimes V_2$ to $V_{m+2}$, the homomorphism $\rho_m^{+}$ is represented as
\begin{equation}
\rho_m^{+}(X)v=\sum_k (v\otimes X, \omega^k_{m+2})v^k_{m+2}, \label{eqn:2-14}
\end{equation}
where $(\cdot,\cdot)$ is inner product on $V_m\otimes V_2$. If we use another orthonormal basis $\{ \rho_{m+2}(g)v^k_{m+2}\}_k$, then we have
\begin{equation}
\rho_m^{+}(X)v=\sum_k (v\otimes X, (\rho_m\otimes \rho_2)(g)\omega^k_{m+2})\rho_{m+2}(g)v^k_{m+2}. \nonumber
\end{equation}
It follows that
\begin{equation}
\begin{split}
\rho_m^+(gXg^{-1})v&=\sum (v\otimes gXg^{-1},\omega^k_{m+2})v^k_{m+2} \\
                   &=\sum (v\otimes \rho_2(g)X,\omega^k_{m+2})v^k_{m+2}  \\
             &=\sum ((\rho_m\otimes \rho_2)(g)(\rho_m(g^{-1})v\otimes X),\omega^k_{m+2})v^k_{m+2} \\
         &=\sum (\rho_m(g^{-1})v\otimes X,(\rho_m\otimes \rho_2)(g^{-1})\omega^k_{m+2})v^k_{m+2} \\
 &=\sum (\rho_m(g^{-1})v\otimes X,\omega^k_{m+2})\rho_{m+2}(g)v^k_{m+2}\\         &=\rho_{m+2}(g)\rho_m^+(X)\rho_m(g^{-1})v.
\end{split} \nonumber
\end{equation}
Thus we have proved the lemma.
\end{proof}
The infinitesimal version of this lemma is given as follows. 
\begin{lem}\label{cor:2-4}
For $X,Y$ in $\mathbf{R}^3\simeq\mathfrak{su}(2)$Cit holds that
\begin{align}
\rho_m^0([X,Y])&=[\rho_m^0(X),\rho_m^0(Y)], \label{eqn:2-15}\\
\rho_m^{\pm}([X,Y])&=\rho_{m\pm 2}^0(X)\rho_{m}^{\pm}(Y)-\rho_m^{\pm}(Y)\rho_m^0(X), \label{eqn:2-16}
\end{align}
where $[\cdot,\cdot]$ denotes the Lie bracket in $\mathfrak{su}(2)$.
\end{lem}

Now, we know that the usual Clifford actions $\{\sigma_i\}=\{e_i\cdot\}$ satisfy the relations 
\begin{equation}
\sigma_i \sigma_j+\sigma_i \sigma_j =-2\delta_{ij}\quad (0\le i, j\le 3). \label{eqn:2-17} 
\end{equation}
We should find what relations the Clifford homomorphisms satisfy. 
\begin{lem}\label{lem:2-5}
The Clifford homomorphisms have the following relations: for $X$, $Y$ in $\mathbf{R}^3\simeq \mathfrak{su}(2)$, 
\begin{align}
\rho_{m+2}^0(X)\rho_{m}^{+}(Y)-\rho_m^{+}(X)\rho_m^0(Y)
    &=\frac{m+2}{2}\rho_m^{+}([X,Y]) ,                  \label{eqn:2-18}\\ 
\rho_{m-2}^0(X)\rho_{m}^{-}(Y)-\rho_m^{-}(X)\rho_m^0(Y)
    &=-\frac{m}{2}\rho_m^{-}([X,Y]),                      \label{eqn:2-19}\\
 \rho_{m}^0(X)\rho_{m}^{0}(Y)+\rho_{m-2}^{+}(X)\rho_m^-(Y)
    &=\frac{m}{2}\rho_m^{0}([X,Y])-m^2(X,Y),                 \label{eqn:2-20} \\
 \rho_{m}^0(X)\rho_{m}^{0}(Y)+\rho_{m+2}^{-}(X)\rho_m^+(Y)
    &=-\frac{m+2}{2}\rho_m^{0}([X,Y])-(m+2)^2(X,Y),           \label{eqn:2-21} 
\end{align}
where $(\cdot,\cdot)$ is the inner product on $\mathbf{R}^3$. 
\end{lem}
\begin{proof}
By direct calculations. 
\end{proof}
We remark that, for $m=1$, the relation \eqref{eqn:2-20} is the usual Clifford relation \eqref{eqn:2-17}.
\section{The higher spin bundles and the higher spin Dirac operators}\label{sec:3}
Let $M$ be the $3$-dimensional oriented Riemannian manifold. Since such a manifold is automatically a spin manifold, we have a principal spin bundle $\mathbf{Spin}(M)$ on $M$, where the structure group is $Spin(3)=SU(2)$. Then all the associated complex vector bundles are induced from the representations of $SU(2)$. For any $m\ge 0$, we define the spin-$\frac{m}{2}$ bundle $\mathbf{S}_m$ by 
\begin{equation}
\mathbf{S}_m=\mathbf{S}_m(M):=\mathbf{Spin}(M)\times_{\rho_m} V_m. \label{eqn:3-1}
\end{equation}
The inner product on $V_m$ induces the one on each fiber of $\mathbf{S}_m$ naturally, which we denote by $\langle\cdot,\cdot\rangle$ on $(\mathbf{S}_m)_x$. 
For example, the spin-$0$ bundle $\mathbf{S}_0$ is the trivial rank $1$ bundle $M\times \mathbf{C}\simeq \Lambda^0(M)\otimes \mathbf{C}$, the spin-$\frac{1}{2}$ bundle $\mathbf{S}_1$ is the spinor bundle, and the spin-$1$ bundle $\mathbf{S}_2$ is $T(M)\otimes \mathbf{C}\simeq \Lambda^1(M)\otimes \mathbf{C}$. 

The spinor bundle $\mathbf{S}_1$ is known as a bundle of modules over the Clifford bundle $\mathbf{C}l(M)$ and the action of $T(M)$ on $\mathbf{S}_1$ is given by 
\begin{equation}
T(M)\times \mathbf{S}_1\ni ([p,e_i],[p,v])\mapsto [p,e_i\cdot v]\in \mathbf{S}_1 , \label{eqn:3-2}
\end{equation}
where $p$ is in $\mathbf{Spin}(M)$, $e_i$ in $\mathbf{R}^3$, and $v$ in $V_1$. 
In the same way, we define the Clifford homomorphisms of $T(M)$ on the higher spin bundle $\mathbf{S}_m$ as follows:
\begin{align}
T(M)\times \mathbf{S}_m\ni ([p,e_i],[p,v])\mapsto [p,\rho_m^0(e_i)v]\in \mathbf{S}_m , \label{eqn:3-3} \\
T(M)\times \mathbf{S}_m\ni ([p,e_i],[p,v])\mapsto [p,\rho_m^{\pm}(e_i)v]\in \mathbf{S}_{m\pm 2} . \label{eqn:3-4}
\end{align}
We can easily check from lemma \ref{lem:2-3} that these bundle homomorphisms are well-defined. 

Before considering the higher spin Dirac operators on $\Gamma (M,\mathbf{S}_m)$, we recall the definition of the Dirac operator $D$ on $\Gamma(M,\mathbf{S}_1)$. Let $\nabla$ be the covariant derivative associated to the spin connection. The Dirac operator $D$ has the following (local) formula:
\begin{equation}
D=\sum_{i=1}^3e_i\cdot\nabla_{e_i}.\label{eqn:3-5}
\end{equation}
On the other hand, we know another description of $D$ as follows: the Dirac operator $D$ is said to be the composed mapping $\mathrm{pr}\circ \nabla$,
\begin{equation}
\Gamma(M,\mathbf{S}_1)\xrightarrow{\nabla}\Gamma(M,\mathbf{S}_1\otimes T^{\ast}(M))\simeq \Gamma(M,\mathbf{S}_1\otimes T(M))\xrightarrow{\mathrm{pr}} \Gamma(M,\mathbf{S}_1),   \label{eqn:3-6}
\end{equation}
where we use $\mathbf{S}_1\otimes T(M)\simeq \mathbf{S}_1\otimes \mathbf{S}_2\simeq \mathbf{S}_3\oplus \mathbf{S}_1$. 

We generalize this composed mapping to give the higher spin Dirac operator (see \cite{BS}, \cite{Bu} and \cite{F}). Since the tensor bundle $\mathbf{S}_m\otimes \mathbf{S}_2$ is isomorphic to $\mathbf{S}_{m+2}\oplus \mathbf{S}_m\oplus \mathbf{S}_{m-2}$, we have three composed mappings for each bundle:
\begin{align}
D^{0}_m:\Gamma(M,\mathbf{S}_m)\xrightarrow{\nabla}\Gamma(M,\mathbf{S}_m\otimes T^{\ast}M) &\xrightarrow{\mathrm{pr}^0} \Gamma(M,\mathbf{S}_m), \label{eqn:3-7} \\
D_m^{\pm}:\Gamma(M,\mathbf{S}_m)\xrightarrow{\nabla}\Gamma(M,\mathbf{S}_m\otimes T^{\ast}M) &\xrightarrow{\mathrm{pr}^{\pm}} \Gamma(M,\mathbf{S}_{m\pm 2}). \label{eqn:3-8}
\end{align}
We call these first order differential operators {\it the higher spin Dirac operators}. In \cite{F}, Fegan show that these operators are conformally invariant first order differential operators and all such operators are given in this way. 

The Clifford homomorphisms in section \ref{sec:2} lead us to represent the higher spin Dirac operators by local formulas such as \eqref{eqn:3-5}.
\begin{prop}\label{prop:3-1}
Let $M$ be the $3$-dimensional spin manifold, $\{ e_i \}_{1\le i \le 3}$ a local orthonormal frame of $T(M)$, and $\nabla$ the covariant derivative associated to the spin connection on $\mathbf{S}_m$. Then we have the following conformally invariant first order differential operators:
\begin{align}
D_m^0&=\sum_{1\le i \le 3}\rho_m^0(e_i)\nabla_{e_i}: \Gamma (M,\mathbf{S}_m)\to \Gamma (M,\mathbf{S}_m),
                           \label{eqn:3-9}   \\
D_m^{\pm}&=\sum_{1\le i \le 3}\rho_m^{\pm}(e_i)\nabla_{e_i}: \Gamma (M,\mathbf{S}_m)\to \Gamma (M,\mathbf{S}_{m\pm 2}).
                              \label{eqn:3-10}
\end{align}
\end{prop}
\begin{ex}
Some higher spin Dirac operators are well-known differential operators. 
\begin{enumerate}
\item 
$D_0^{+}$ is $2d$ on $\Gamma(M,\mathbf{S}_0)=\Gamma(M,\Lambda^0(M)\otimes \mathbf{C})$
 \item 
$D^0_1$ is the Dirac operator $D$ and $D_1^+$ is the twistor operator on $\Gamma(M,\mathbf{S}_1)$. 
\item $D_2^0$ is $2\ast d$ and  $D_2^-$ is $2d^{\ast}$ on $\Gamma(M,\mathbf{S}_2)=\Gamma(M,\Lambda^1(M)\otimes \mathbf{C})$, where $\ast$ is the Hodge star operator from $\Lambda^1(M)$ to $\Lambda^2(M)$. 
\end{enumerate}
\end{ex}
From the discussion in section \ref{sec:2}, we can derive some properties of the higher spin Dirac operators. 

First, we discuss the adjointness of the operators. On $\Gamma(M,\mathbf{S}_m)$, we set the inner product by
\begin{equation}
(\phi_1,\phi_2):=\int_M\langle \phi_1(x),\phi_2(x)\rangle dx,
\end{equation}
where $dx$ denotes the volume element of $M$. 
\begin{prop}\label{lem:3-2}
We denote the formal adjoint of a differential operator $A$ by $A^{\ast}$. Then we have 
\begin{align}
(D_m^0)^{\ast}&=D_m^0,                       \label{eqn:3-11}\\
(D_m^{\pm})^{\ast}&=D_{m\pm 2}^{\mp}.        \label{eqn:3-12}
\end{align}
In particular, $D_m^0$ is formally self adjoint.
\end{prop}
\begin{proof}
We can easily show that the Dirac operator is formally self-adjoint (for example, see \cite{LM}). In the same way, we can prove \eqref{eqn:3-11} and \eqref{eqn:3-12} by using lemma \ref{lem:2-2},.
\end{proof}
Next, we shall discuss the commutativity among the operators. So we have to introduce some curvature homomorphisms. For vector fields $X$, $Y$, the curvature $R_m$ for $\mathbf{S}_m$ is given by
\begin{align}
R_m(X,Y)&=\nabla_X \nabla_Y-\nabla_Y \nabla_X-\nabla_{[X,Y]} \quad \in \Gamma (M,\mathrm{End}(\mathbf{S}_m))\label{eqn:3-13} \\
&=\frac{1}{4}\sum_{\sigma \in S_3}\mathrm{sgn}(\sigma)\langle R(X,Y)(e_{\sigma (1)}),e_{\sigma (2)}\rangle\rho_m^0(e_{\sigma (3)}), \label{eqn:3-14}
\end{align}
where $R(\cdot ,\cdot)$ is the curvature transformation for $T(M)$ and $\{e_i\}_{1\le i \le 3}$ is a local orthonormal frame on $T(M)$. Then we obtain the following curvature homomorphisms from $\mathbf{S}_m$ to $\mathbf{S}_m$ or $\mathbf{S}_{m\pm 2}$:
\begin{align}
R_m^0:&=\sum_{\sigma \in S_3}\mathrm{sgn}(\sigma)
             \rho_m^0(e_{\sigma (1)})R_m(e_{\sigma (2)},e_{\sigma (3)})
                   \in \Gamma (M,\mathrm{End}(\mathbf{S}_m)), \label{eqn:3-15}\\
R_m^{\pm}:&=\sum_{\sigma \in S_3}\mathrm{sgn}(\sigma)
             \rho_m^{\pm}(e_{\sigma (1)})R_m(e_{\sigma (2)},e_{\sigma (3)})
                 \in \Gamma (M,\mathrm{Hom}(\mathbf{S}_m,\mathbf{S}_{m\pm 2})) .
                                                     \label{eqn:3-16}
\end{align}
Here, we show that $(R_m^0)^{\ast}=R_m^0$ and $(R_m^{\pm})^{\ast}=R_{m\pm 2}^{\mp}$. In particular, $(R_m^0)_x$ has real eigenvalues for each $x$ in $M$.
\begin{ex}
Let $\mathrm{Ric}$ be the Ricci curvature and $\kappa$ the scalar curvature. Then we have
\begin{equation}
R_1^0=\frac{1}{2}\kappa , \quad R_2^0=4\mathrm{Ric}, \quad R_0^0=R_2^-=R_0^+=0 .
                     \label{eqn:3-17}
\end{equation}
\end{ex}
The commutativity among the higher spin Dirac operators follows from Lemma \ref{lem:2-5}. The important fact is that we have two Laplace type operators on $\Gamma(M,\mathbf{S}_m)$ for each $m\ge 1$. 
\begin{thm}\label{thm:3-3}
Let $\nabla^{\ast}\nabla$ be the connection Laplacian on $\mathbf{S}_m$. Then the higher spin Dirac operators satisfy the following Bochner type identities:
\begin{align}
D_m^0D_m^0+D_{m-2}^+D_m^-&=m^2\nabla^{\ast}\nabla+\frac{m}{2}R_m^0,
                    \label{eqn:3-18}\\
D_m^0D_m^0+D_{m+2}^-D_m^+&=(m+2)^2\nabla^{\ast}\nabla-\frac{m+2}{2}R_m^0,
                \label{eqn:3-19} \\
D_{m+2}^0D_m^+-D_{m}^+D_m^0&=\frac{m+2}{2}R_m^+, \label{eqn:3-20} \\
D_{m-2}^0D_m^--D_{m}^-D_m^0&=-\frac{m}{2}R_m^-, \label{eqn:3-21}
\end{align}
\end{thm}
\begin{proof}
We shall prove \eqref{eqn:3-18}. We fix $x$ in $M$ and choose an orthonormal frame $\{e_i\}$ in a neighborhood of $x$ such that $(\nabla_{e_i}e_j)_x=0$ for all $i,j$. Hence, we have $(\nabla_{e_i}\rho_m^0(e_j))_x=0$ for all $i,j$. Then it holds from Lemma \ref{lem:2-5} that
\begin{equation}
\begin{split}
D_m^0D_m^0&+D_{m-2}^+D_m^- \\
&=\sum_{i,j}\left( \rho_m^0(e_i)\nabla_{e_i}\rho_m^0(e_j)\nabla_{e_j}+\rho_{m-2}^+(e_i)\nabla_{e_i}\rho_m^-(e_j)\nabla_{e_j} \right) \\
&=\sum_{i}\left( \rho_m^0(e_i)\rho_m^0(e_i)+\rho_{m-2}^+(e_i)\rho_m^-(e_i) \right)\nabla_{e_i}\nabla_{e_i} \\
& \quad \quad \sum_{i\neq j}\left( \rho_m^0(e_i)\rho_m^0(e_j)+\rho_{m-2}^+(e_i)\rho_m^-(e_j) \right)\nabla_{e_i}\nabla_{e_j} \\
&=-m^2\sum_i \nabla_{e_i}\nabla_{e_i}+\frac{m}{2}\sum_{i<j}\rho_m^0(e_i e_j-e_j e_i)( \nabla_{e_i}\nabla_{e_j}- \nabla_{e_j}\nabla_{e_i}) \\
&=m^2 \nabla^{\ast} \nabla+\frac{m}{2}R_m^0.
 \end{split} \nonumber 
\end{equation}
\end{proof}
\begin{ex}
\begin{enumerate}
\item (the case of $m=0$) The relation \eqref{eqn:3-19} means $d^{\ast}d=\nabla^{\ast}\nabla$ and the relation \eqref{eqn:3-20} does $dd=0$.
\item (the case of $m=1$) The relation \eqref{eqn:3-18} means 
\begin{equation}
D^2=\nabla^{\ast} \nabla+\frac{1}{4}\kappa. \label{eqn:3-22}
\end{equation}
\item (the case of $m=2$) The relation \eqref{eqn:3-18} means 
\begin{equation}
d^{\ast}d+dd^{\ast}=\nabla^{\ast}\nabla +\mathrm{Ric}. \label{eqn:3-23}
\end{equation}
and \eqref{eqn:3-21} does $dd=0$. 
\end{enumerate}
\end{ex}
Now, we denote the Laplace type operators in \eqref{eqn:3-18} and \eqref{eqn:3-19} by 
\begin{align}
\Delta_m:&=D_m^0D_m^0+D_{m-2}^+D_m^{-}, \label{eqn:3-24}  \\
\widetilde{\Delta}_m:&=D_m^0D_m^0+D_{m+2}^-D_m^{+}.\label{eqn:3-25} 
\end{align}
If $M$ is compact, these Laplace type operators are non-negative operators and 
satisfy that 
\begin{gather}
\ker \Delta_m =\ker D_m^0\cap \ker D_m^-,\label{eqn:3-26}\\
\ker \widetilde{\Delta}_m =\ker D_m^0\cap \ker D_m^+,\label{eqn:3-27} \\
\ker \nabla =\ker \Delta_m \cap \ker \widetilde{\Delta}_m=\ker D_m^0\cap \ker D_m^-\cap \ker D_m^+. \label{eqn:3-28}
\end{gather}
The following corollary is the key to give lower bounds for the first eigenvalues of $\Delta_m$ and $\widetilde{\Delta}_m$.
\begin{cor}\label{cor:3-4}
The Laplace type operators $\Delta_m$ and $\widetilde{\Delta}_m$ satisfy that
\begin{equation}
(m+2)^2\Delta_m-m^2\widetilde{\Delta}_m=m(m+1)(m+2)R_m^0  \label{eqn:3-29}
\end{equation}
\end{cor}
\begin{proof}
We eliminate the connection Laplacian $\nabla^{\ast}\nabla$ from \eqref{eqn:3-18} and \eqref{eqn:3-19}.
\end{proof}

Finally, we discuss the ellipticity of the operators. Of course, it is clear that $\Delta_m$ and $\widetilde{\Delta}_m$ are elliptic. 
\begin{prop}\label{prop:3-5}
\begin{enumerate}
\item
The second order differential operator $D_{m+2}^-D_m^+=(D_m^+)^{\ast}D_m^+$ is elliptic for each $m$.\item
If $m$ is odd, then the first order differential operator $D_m^0$ is elliptic. Hence $D_m^0$ is an elliptic self adjoint operator. 
\end{enumerate}
\end{prop}
\begin{proof}
We investigate the ellipticity of $D_m^0$. The principal symbol of $D_m^0$ is
\begin{equation}
\sigma_{\xi} (D_m^0)=\rho_m^0(\xi) ,  \label{eqn:3-30}
\end{equation}
where $\xi=\sum \xi_i e_i$ is in $T_x^{\ast}(M)\simeq T_x(M)$. There exists $g$ in $SU(2)$ such that
\begin{equation}
g\xi g^{-1}=(\xi_1^2+\xi_2^2+\xi_3^2)^{\frac{1}{2}}e_1.  \label{eqn:3-31}
\end{equation}
Then we have 
\begin{equation}
\begin{split}
\det \sigma_{\xi} (D_m^0) &=\det \rho_m(g)\rho_m^0(\xi)\rho_m(g^{-1}) \\
     &=\det \rho_m^0(g\xi g^{-1}) \\
    &=\det \rho_m^0((\xi_1^2+\xi_2^2+\xi_3^2)^{\frac{1}{2}}e_1) \\
    &=(\xi_1^2+\xi_2^2+\xi_3^2)^{\frac{m+1}{2}}\det \rho_m^0(e_1)  \\
    &=(\xi_1^2+\xi_2^2+\xi_3^2)^{\frac{m+1}{2}}\prod_{k=0}^mi(2k-m). 
\end{split}  \label{eqn:3-32}
\end{equation}
It follows that, if $m$ is odd, then $\det \sigma_{\xi} (D_m^0)$ is not zero for $\xi\neq 0$. Hence $D_m^0$ is elliptic. In the same way, we verify that $D_{m+2}^-D_m^+$ is elliptic.
\end{proof}
\begin{cor}
We assume that the spin manifold $M$ is compact. Then $\ker D_m^+$ and $\ker D_{2p+1}$ are finite dimensional vector spaces for any $m$ and $p$. 
\end{cor}
\section{lower bounds for the first eigenvalues of the higher spin Dirac operators}\label{sec:4}
In this section, we assume that $M$ is a $3$-dimensional connected compact spin manifold
. From corollary \ref{cor:3-4},
we have 
\begin{equation}
(m+2)^2(\Delta_m\phi,\phi)-m^2(\widetilde{\Delta}_m\phi,\phi)=m(m+1)(m+2)(R_m^0(\phi),\phi), \label{eqn:4-1}
\end{equation}
where $\phi$ is a section of $\mathbf{S}_m$ and 
\begin{equation}
(R_m^0(\phi),\phi):=\int_M \langle(R_m^0)_x\phi(x),\phi(x)\rangle dx.    \label{eqn:4-2}
\end{equation}
From the above equation \eqref{eqn:4-1}, we can obtain lower bounds estimations for the eigenvalues of $\Delta_m$ and $\widetilde{\Delta}_m$ depending on the curvature transformation $R_m^0$.

First, we consider a lower bound for the first eigenvalue of the Dirac operator $D=D_1^0$. It follows from \eqref{eqn:4-1} that, for a spinor $\phi$ in $\Gamma(M,\mathbf{S}_1)$, 
\begin{align}
& 9\| D\phi \|^2-(\| D\phi \|^2+\| D_1^{+}\phi \|^2)\nonumber \\
=& 8\| D\phi \|^2 -\| D_1^{+}\phi \|^2 \label{eqn:4-3} \\
=& 6(R_1^0(\phi),\phi)=3(\kappa \phi,\phi). \nonumber
\end{align}
Because of $\| D_1^{+}\phi \| \ge 0$, we have 
\begin{equation}
\| D\phi \|^2\ge \frac{3}{8}(\kappa \phi,\phi).  \label{eqn:4-4}
\end{equation}
If $\phi_1$ is an eigenspinor with the first eigenvalue $\lambda_1$ of $D$, 
then $(\lambda_1)^2$ has a lower bound,
\begin{equation}
(\lambda_1)^2 \ge \frac{3(\kappa \phi_1,\phi_1)}{8\|\phi_1\|^2}\ge \frac{3}{8}\kappa_-. 
                   \label{eqn:4-5}
\end{equation}
where
\begin{equation}
\kappa_-:=\min_{x \in M} \kappa (x).   \label{eqn:4-6}
\end{equation}
If the equality holds in \eqref{eqn:4-5}, then $\phi_1$ is in $\ker D_1^{+}$, that is, $\phi_1$ is a twistor spinor. This inequality coincides with the ones given by Friedrich (see \cite{BF}). 

Next, we investigate the case of the elliptic operator $D_3^-D_1^+=(D_1^+)^{\ast}D_1^+$. It holds that 
\begin{equation}
\begin{split}
(D_3^-D_1^+\phi,\phi) &=8\|D\phi\|^2-3(\kappa \phi, \phi) \\
                 &\ge -3(\kappa \phi ,\phi).  \label{eqn:4-7}
      \end{split}
   \end{equation}
If we denote the first eigenvalue of $D_3^-D_1^+$ by $\mu_1$, then we have 
\begin{equation}
\mu_1 \ge -3\kappa_+,       \label{eqn:4-8}
\end{equation}
where
\begin{equation}
\kappa_+=\max_{x \in M}\kappa (x). \label{eqn:4-9}
\end{equation}

In general case ($m\ge 2$), we have the inequalities 
\begin{gather}
( \Delta_m\phi, \phi)=\|D_m^0\phi \|^2+\|D_m^-\phi \|^2
\ge \frac{m(m+1)}{m+2}(\mathrm{R}_m^0(\phi),\phi ),    \label{eqn:4-10} \\
( \widetilde{\Delta}_m\phi, \phi)=\left(\|D_m^0\phi \|^2+\|D_m^+\phi \|^2\right)\ge -\frac{(m+2)(m+1)}{m}(\mathrm{R}_m^0(\phi),\phi ) .  \label{eqn:4-11}
\end{gather}
Then we give lower bounds for the first eigenvalues of $\Delta_m$ and $\widetilde{\Delta}_m$.
\begin{thm}\label{thm:4-1} We assume that there exist constants $r_{m-}$ and $r_{m+}$ such that 
\begin{equation}
r_{m-}\|\phi\|^2\le (R_m^0(\phi),\phi)\le r_{m+}\|\phi\|^2
\end{equation} 
for any $\phi$ in $\Gamma(M,\mathbf{S}_m)$.
\begin{enumerate}
\item Let $\lambda_1$ be the first eigenvalue of $\Delta_m$. Then we have the inequality
\begin{equation}
\lambda_1 \ge \frac{m(m+1)}{m+2}  r_{m-}.
                             \label{eqn:4-12}
\end{equation}
If the equality holds in \eqref{eqn:4-12}, the eigenvectors with the eigenvalue $\lambda_1$  is in $\ker \widetilde{\Delta}_m$. 
\item 
Let $\mu_1$ be the first eigenvalue of $\widetilde{\Delta}_m$. Then we have the inequality 
\begin{equation}
\mu_1 \ge -\frac{(m+2)(m+1)}{m}  r_{m+}. \label{eqn:4-14}                           
\end{equation}
If the equality holds in \eqref{eqn:4-14}, then the eigenvectors with the eigenvalue $\mu_1$ is in $\ker \Delta_m$. 
\end{enumerate}
\end{thm}
\begin{cor}[\cite{GM}]\label{cor:4-2}
We assume that there exists a constant $\mathrm{ric}_-$ such that 
\begin{equation}
(\mathrm{Ric}(\phi),\phi)\ge \mathrm{ric}_-\|\phi\|^2
\end{equation}
for any $\phi$ in $\Gamma(M,\Lambda^1(M))$. 
Let $\lambda_1$ be the first eigenvalue of the Laplace-Beltrami operator $dd^{\ast}+d^{\ast}d$ on $\Gamma(M,\Lambda^1(M))$. Then we have 
\begin{equation}
\lambda_1 \ge \frac{3}{2}\mathrm{ric}_-.  \label{eqn:4-16}
\end{equation}
If the equality holds in \eqref{eqn:4-16}, the eigenforms with the eigenvalue $\lambda_1$ is in $\ker \widetilde{\Delta}_2=\ker d \cap \ker D_2^+$.
\end{cor}
\section{On the $3$-dimensional manifold of constant curvature}\label{sec:5}
In this section, we shall discuss the higher spin Dirac operators on $3$-dimensional manifold of constant curvature. 
\begin{lem}\label{lem:5-1}
On the $3$-dimensional spin manifold $M$ of constant curvature $c$, the curvature homomorphism $R^0_m$ is $m(m+2)c$ and $R_m^{\pm}$ is zero.
\end{lem}
\begin{proof}
Since $M$ has constant curvature, it holds that, for vector fields $X$, $Y$, and $Z$,
$$R(X,Y)Z=c\{ (Y,Z)X-(X,Z)Y\}.$$
Then we have $(R(e_i,e_j)e_k,e_l)=c(\delta_{j k}\delta_{i l}-\delta_{i k}\delta_{j l})$. Hence, 
$$R^0_m=-c\sum \rho_m^0(e_i)\rho_m^0(e_i)=m(m+2)c, \quad R_m^{\pm}=0.$$ 
Here, we use that $-\sum \rho_m^0(\sigma_i)\rho_m^0(\sigma_i)$ is the Casimir operator on $V_m$. 
\end{proof}
\begin{prop}\label{prop:5-2}
On the $3$-dimensional spin manifold of constant curvature $c$, it hold that 
\begin{align}
D_m^0D_m^0+D_{m-2}^+D_m^-&=m^2\nabla^{\ast}\nabla+\frac{m^2(m+2)}{2}c,
                    \label{eqn:7-2}\\
D_m^0D_m^0+D_{m+2}^-D_m^+&=(m+2)^2\nabla^{\ast}\nabla-\frac{m(m+2)^2}{2}c,
                \label{eqn:7-3} \\
D_{m+2}^0D_m^+-D_{m}^+D_m^0&=0,\label{eqn:7-4} \\
D_{m-2}^0D_m^--D_{m}^-D_m^0&=0. \label{eqn:7-5}
\end{align}
In particular, we have
\begin{gather}
\Delta_m D_m^0=D_m^0\Delta_m, \quad \widetilde{\Delta}_mD_m^0=D_m^0\widetilde{\Delta}_m, \\
\Delta_m(D_{m-2}^+D_m^-)=(D_{m-2}^+D_m^-)\Delta_m, \quad \widetilde{\Delta}_m(D_{m+2}^-D_m^+)=(D_{m+2}^-D_m^+)\widetilde{\Delta}_m.
\end{gather}
\end{prop}
We conclude from this proposition that $\Delta_m$, $D_m^0$, and $D_{m-2}^+D_m^-$ have the simultaneous eigenspaces. As an example, we will calculate the eigenvalues of these operators on $S^3$ in the next section.
\section{The spectra of the higher spin Dirac operators on $S^3$}\label{sec:6}
In this section, we calculate all the eigenvalues of the higher spin Dirac operators on the symmetric space $S^3$ with constant curvature $1$. In \cite{H}, the author gives a method for calculating of the eigenvalues and the eigenspinors for the Dirac operator on $S^3$. We can use the method in our situation and calculate the eigenvalues. So we refer to the paper \cite{H} for details. 

First, we shall explain the $3$-dimensional sphere $S^3$ as the symmetric space $Spin(4)/Spin(3)$. It is well-known that $Spin(4)$ and $Spin(3)$ are isomorphic to $SU(2)\times SU(2)$ and $SU(2)$, respectively. We realize $S^3$ as $SU(2)$, 
\begin{equation}
    S^3\ni x=(x_1,x_2,x_3,x_4)\mapsto h=
            \begin{pmatrix}
                x_4+ix_1 & x_2+ix_3 \\
                     -x_2+ix_3 & x_4-ix_1 
        \end{pmatrix}     \label{eqn:6-1}
     \in SU(2).     
\end{equation}
Therefore, the action of $SU(2)\times SU(2)$ on $S^3$ is represented by 
\begin{equation}
    (SU(2)\times SU(2))\times S^3\ni(g,h)\mapsto phq^{-1}\in S^3, \label{eqn:6-2}
             \end{equation}
where $g=(p,q)$ is in $SU(2)\times SU(2)$. 
Since the isotropy subgroup of $e=(0,0,0,1)$ is the subgroup $SU(2)$ in $SU(2)\times SU(2)$, we have the symmetric space $S^3$,
\begin{equation}
S^3=Spin(4)/Spin(3)=SU(2)\times SU(2)/\mathrm{diag} SU(2). \label{eqn:6-3}
\end{equation}
Here, the map `$\mathrm{diag}$' is given by
\begin{equation}
\mathrm{diag}: SU(2)\ni h \mapsto (h,h)\in SU(2)\times SU(2).  \label{eqn:6-4}
\end{equation}
Besides, the principal spin bundle $\mathbf{Spin}(S^3)$ is the Lie group $Spin(4)$, whose projection from the total space to the base space is
\begin{equation}
\mathbf{Spin}(S^3)=Spin(4)\ni g\mapsto pq^{-1}\in S^3.  \label{eqn:6-5}
\end{equation}
This principal spin bundle induces the spin $\frac{m}{2}$ bundle $\mathbf{S}_m$ as a homogeneous vector bundle: 
\begin{equation}
\mathbf{S}_m:=Spin(4)\times_{\rho_m} V_m. \label{eqn:6-6}
\end{equation}
Hence the space of sections $L^2(S^3, \mathbf{S}_m)$ is a representation space of $Spin(4)$.

Now, we trivialize the vector bundle $\mathbf{S}_m$ as follows:
\begin{equation}
   \mathbf{S}_m=Spin(4)\times_{\rho_m}
               V_m \ni [g, v] \mapsto (pq^{-1},
                         \rho_m(p)v)\in S^3\times V_m. 
                            \label{eqn:6-7}
   \end{equation}
So the sections of $\mathbf{S}_m$ are represented as the $\mathbf{C}^{m+1}$-valued or the $V_m$-valued functions on $S^3$. 
In this situation, we can present explicit formulas of the higher spin Dirac operators on $S^3$, where the operators acts on the $V_m$-valued functions.
\begin{prop}\label{prop:6-1}
For the trivialization \eqref{eqn:6-7}, the higher spin Dirac operators on $S^3$ are the following:
\begin{align}
D_m^0&=\frac{m(m+2)}{2}+\sum \rho_{m}^0(e_i)Z_i, \label{eqn:6-8} \\
D^{\pm}_m&=\sum \rho_{m}^{\pm}(e_i)Z_i.  \label{eqn:6-9}
\end{align}
Here $Z_i$ is the right invariant vector field on $S^3=SU(2)$ corresponding to $\sigma_i$ in $\mathfrak{su}(2)$, which is given by  
\begin{gather}
  Z_1 = -x_1\frac{\partial}{\partial x_4}
             +x_4\frac{\partial}{\partial x_1}
                  -x_3\frac{\partial}{\partial x_2}
                     +x_2\frac{\partial}{\partial x_3},\nonumber \\
     Z_2 = -x_2\frac{\partial}{\partial x_4}
               +x_3\frac{\partial}{\partial x_1}
                    +x_4\frac{\partial}{\partial x_2}
                          -x_1\frac{\partial}{\partial x_3}, 
   \label{eqn:6-10}  \\
   Z_3 = -x_3\frac{\partial}{\partial x_4}
            -x_2\frac{\partial}{\partial x_1}
                  +x_1\frac{\partial}{\partial x_2}
                      +x_4\frac{\partial}{\partial x_3}.\nonumber
\end{gather}
\end{prop} 
\begin{cor}\label{cor:6-2}
The Laplace type operators $\Delta_m$ and $\widetilde{\Delta}_m$ are represented by
\begin{align}
\Delta_m &= -m^2\sum Z_i^2+m^2D_m^0-\frac{m^2(m+2)(m-2)}{4},  \label{eqn:6-11}\\
\widetilde{\Delta}_m &= -(m+2)^2\sum Z_i^2+(m+2)^2D_m^0-\frac{m(m+2)^2(m+4)}{4}. \label{eqn:6-12}
\end{align}
\end{cor}

Since the higher spin Dirac operators on $S^3$ are homogeneous differential operators, the eigenspaces are representation spaces of $Spin(4)$. So we have to decompose $L^2(S^3,\mathbf{S}_m)$ into its irreducible components. By the Frobenius reciprocity, we have the following lemma.
\begin{lem}\label{lem:6-3}
The representation space $L^2(S^3,\mathbf{S}_{m})$ decomposes into its irreducible components as follows: 
\begin{enumerate}
\item (the case of $m=2p+1$)
\begin{equation}
L^2(S^3,\mathbf{S}_{2p+1})  \simeq \bigoplus_{0\le s \le p \atop k \ge p-s} E_{k,k+2s+1}\oplus E_{k+2s+1, k}. \label{eqn:6-13}
\end{equation}
\item  (the case of $m=2p$) 
\begin{equation}
L^2(S^3,\mathbf{S}_{2p})  \simeq \bigoplus_{1\le s \le p \atop k \ge p-s} E_{k,k+2s}\oplus E_{k+2s, k}\bigoplus_{k\ge p}E_{k,k}.\label{eqn:6-14}
\end{equation}
\end{enumerate}
Here $E_{k,l}$ is the representation space for the exterior tensor product representation $\rho_k\widehat{\otimes}\rho_l$ of $Spin(4)=SU(2)\times SU(2)$ and $\dim E_{k,l}=(k+1)(l+1)$.  
\end{lem}
We calculate the action of the higher spin Dirac operators on $E_{k,l}$ by the method given in \cite{H}. Then we have the following propositions.
\begin{prop}\label{thm:6-4}
\begin{enumerate}
\item
The eigenvalues of the self adjoint operator $D_{m}^0$ on $S^3$ are given as follows: \begin{enumerate}
\item (the case of $m=2p+1$) 
\begin{equation}
\begin{cases}
(2s+1)(k+\frac{2s+3}{2}) &  \mbox{on $E_{k+2s+1,k}$}, \\
-(2s+1)(k+\frac{2s+3}{2}) & \mbox{on $E_{k,k+2s+1}$}. 
\end{cases}\label{eqn:6-15}
\end{equation}
In particular, $\ker D_{2p+1}^0$ is zero for each $p$.
\item (the case of $m=2p$) 
\begin{equation}
\begin{cases}
2s(k+s+1) &  \mbox{on $E_{k+2s,k}$}, \\
-2s(k+s+1) & \mbox{on $E_{k,k+2s}$}, \\
0 & \mbox{on $E_{k,k}$}.
\end{cases}\label{eqn:6-16}
\end{equation}
\end{enumerate}
\item
The eigenvalues of the second order operator $D^{+}_{m-2}D_m^-$ on $S^3$ are given as follows: 
\begin{enumerate}
\item (the case of $m=2p+1$) 
\begin{multline}
4(p-s)(k+1-(p-s))(p+s+1)(k+p+s+2)\\
 \mbox{on $E_{k+2s+1,k}$ or $E_{k,k+2s+1}$}.
\end{multline}
\item (the case of $m=2p$) 
\begin{multline}
4(p-s)(k+1-(p-s))(p+s)(k+p+s+1)\\
 \mbox{on $E_{k+2s,k}$ or $E_{k,k+2s}$}. 
\end{multline}
\end{enumerate}
\item
The eigenvalues of the second order elliptic operator $D^{-}_{m+2}D_m^+$ on $S^3$ are given as follows: 
\begin{enumerate}
\item (the case of $m=2p+1$) 
\begin{multline}
4(p-s+1)(k-(p-s))(p+s+2)(k+p+s+3) \\
    \quad \mbox{on $E_{k+2s+1,k}$ or $E_{k,k+2s+1}$}.
\end{multline}
\item (the case of $m=2p$) 
\begin{multline}
4(p-s+1)(k-(p-s))(p+s+1)(k+p+s+2)\\
  \mbox{on $E_{k+2s,k}$ or $E_{k,k+2s}$}.
\end{multline}
\end{enumerate}
In particular,
\begin{equation}
\dim \ker D_{m}^+=\frac{1}{6}(m+1)(m+2)(m+3). 
\end{equation}
\end{enumerate}
\end{prop}
\begin{prop}
\begin{enumerate}
\item
The eigenvalues of the Laplace type operator $\Delta_m$ are as follows: 
\begin{enumerate}
\item (the case of $m=2p+1$) 
\begin{multline}
(2s+1)^2(k+\frac{2s+3}{2})^2\\
+4(p-s)(k+1-(p-s))(p+s+1)(k+p+s+2)\\
 \quad \mbox{on $E_{k+2s+1,k}$ or $E_{k,k+2s+1}$}.
\end{multline}
\item (the case of $m=2p$)
\begin{multline}
(2s)^2(k+s+1)^2\\
 +4(p-s)(k+1-(p-s))(p+s)(k+p+s+1)\\
  \mbox{on $E_{k+2s,k}$ or $E_{k,k+2s}$}.
\end{multline}
\end{enumerate}
In particular, $\ker \Delta_{m}$ is zero. 
\item The eigenvalues of the Laplace type operator $\widetilde{\Delta}_m$ are given as follows:
\begin{enumerate}
\item (the case of $m=2p+1$) 
\begin{multline}
(2s+1)^2(k+\frac{2s+3}{2})^2 \\
  +4(p-s+1)(k-(p-s))(p+s+2)(k+p+s+3) \\
 \quad \mbox{on $E_{k+2s+1,k}$ or $E_{k,k+2s+1}$}.
 \end{multline}
In particular, $\ker \widetilde{\Delta}_{2p+1}$ is zero. 
\item (the case of $m=2p$)
\begin{multline}
(2s)^2(k+s+1)^2  \\
    +4(p-s+1)(k-(p-s))(p+s+1)(k+p+s+2)   \\
 \quad  \mbox{on $E_{k+2s,k}$ or $E_{k,k+2s}$}.
\end{multline}
In particular, $\dim \ker \widetilde{\Delta}_{2p}=(p+1)^2$.
\end{enumerate}
\end{enumerate}
\end{prop}
\begin{rem}
For the $3$-dimensional flat manifold, that is, $T^3=S^1\times S^1\times S^1$, we can easily calculate the dimensions of the kernels for the higher spin Dirac operators:
\begin{gather}
\dim \ker D_m^+=\dim \ker \Delta_m=\dim \ker \widetilde{\Delta}_m=m+1, \\
\dim \ker D_{2p+1}^0=2p+2.
\end{gather}
\end{rem}
\section*{Acknowledgements}
The author would like to thank T. Kori for his encouragement.

\end{document}